\documentclass[12pt,reqno, oneside]{amsart}
\usepackage{graphicx}

\def\comment#1{(\textbf{Comment: }\textsl{#1})}
\def\comment#1{}

\newtheorem{theorem}{Theorem}[section]
\newtheorem{cor}[theorem]{Corollary}
\newtheorem{lemma}[theorem]{Lemma}

\numberwithin{equation}{section}

\def\dpow#1#2{\frac{#1^{#2}}{#2!}}

\def\sumz#1{\sum_{#1=0}^\infty}
\def\sumo#1{\sum_{#1=1}^\infty}
\def\egf#1{\sumz n #1\dpow xn}

\def\dpath{decomposing path}
\DeclareMathOperator{\prop}{prop}
\DeclareMathOperator{\tree}{tree}
\DeclareMathOperator{\lucky}{lucky}
\DeclareMathOperator{\asc}{asc}
\DeclareMathOperator{\des}{des}
\DeclareMathOperator{\comp}{comp}

\newcommand{\newsubsection}{\medskip}

\begin{document}
\title[A refinement of Cayley's formula]{A refinement of Cayley's formula for trees}
\author{Ira M. Gessel$^*$ and Seunghyun Seo}
\address{Department of Mathematics\\
   Brandeis University\\
   Waltham, MA 02454-9110}
\email{gessel@brandeis.edu}
\address{Department of Mathematics\\
Korea Advanced Institute of Science and Technology\\
Daejon 305-701\\
Korea}
\email{shseo@kaist.ac.kr}
\date{July 17, 2005}
\thanks{$^*$Partially supported by NSF Grant DMS-0200596}
\begin{abstract}
A \emph{proper vertex} of a rooted tree with totally ordered vertices is a vertex that is less than all its proper descendants. We count several kinds of labeled rooted trees and forests  by the number of proper vertices. Our results are all  expressed in terms of the polynomials 
\[P_n(a,b,c)= c\prod_{i=1}^{n-1}(ia+(n-i)b +c),\]
which reduce to to $(n+1)^{n-1}$ for $a=b=c=1$.

Our study of proper vertices was motivated by 
Postnikov's hook length formula
\[(n+1)^{n-1}=\frac{n!}{2^n}\sum _T \prod_{v}\left(1+\frac1{h(v)}\right),\]
where the sum is over all unlabeled binary trees $T$ on $n$ vertices, the product is over all vertices $v$ of $T$, and $h(v)$ is the number of descendants of $v$. Our results give analogues of Postnikov's formula for other types of trees, and we also find an interpretation of the polynomials $P_n(a,b,c)$ in terms of parking functions.
\end{abstract}

\maketitle
\section{Introduction}

Cayley \cite{cayley} showed that there are $n^{n-2}$ (unrooted) trees on $n$ vertices. Equivalently, 
there are $(n+1)^{n-1}$ forests of rooted labeled trees on $n$ vertices. In this paper 
we study the homogeneous polynomials
\begin{equation}
\label{e-0}
P_n(a,b,c)= c\prod_{i=1}^{n-1}(ia+(n-i)b +c).
\end{equation}
which reduce to $(n+1)^{n-1}$ for $a=b=c=1$. We refine Cayley's formula by showing that $P_n(a,b,c)$ counts rooted forests by the number of trees and the number of \emph{proper vertices}, which are vertices 
that are less than all of their proper descendants. Moreover, other
evaluations of $P_n(a,b,c)$ have similar interpretations for other types of trees and forests: $k$-ary trees, forests of ordered trees, and forests of $k$-colored ordered trees (which reduce to ordered trees for $k=1$). We also give an interpretation of $P_n(a,b,c)$ in terms of parking functions.



This work developed from a problem posed by Alexander Postnikov at the Richard Stanley 60th Birthday Conference. He gave the following identity, which he had derived indirectly \cite{post1, post2} and asked for a direct proof:
\begin{equation}
\label{e-P1}
(n+1)^{n-1}=\frac{n!}{2^n}\sum _T \prod_{v}\left(1+\frac1{h(v)}\right),
\end{equation}
where the sum is over all unlabeled (incomplete) binary trees $T$ on $n$ vertices, the product is over all vertices $v$ of $T$, and $h(v)$ is the \emph{hook length} of $v$, which is the number of descendants of $v$, including $v$.

Direct combinatorial proofs of Postnikov's identity were given by 
Seo \cite{seo}, Du and Liu \cite{du-liu}, and Chen and Yang \cite{chen-yang}. Seo's approach involved  counting binary trees by proper vertices in the following way:
We rewrite \eqref{e-P1} as
\begin{equation}
\label{e-P2}
2^n(n+1)^{n-1}=n!\sum _T \prod_{v}\left(1+\frac1{h(v)}\right),
\end{equation}

As we will see in Section \ref{s-hook} (see also \cite{seo}), a straightforward argument shows that 
the right side of \eqref{e-P2} is equal to
$
\sum_B 2^{\prop B}
$ 
where the sum is over all labeled binary trees $B$ on $[n]=\{1,2,\dots,n\}$ and $\prop B$ is the number of proper vertices of $B$. This observation led us to consider the polynomial 
$\sum_B t^{\prop B}$, which was shown in \cite{seo} to be equal to 
$P_n(2,t,t).$ Extending this result to other kinds of forests and trees by proper vertices, and we found, surprisingly, that all of the results can be expressed in terms of the polynomials 
$P_n(a,b,c)$. Moreover, each type of forest or tree has a hook length formula analogous to Postnikov's. (Some of these hook length formulas have been found earlier by Du and Liu \cite{du-liu}.)

A closely related statistic to the number of proper vertices is the number of \emph{proper edges}, studied by Shor \cite{shor},  Zeng \cite{zeng}, and Chen and Guo \cite{chen-guo}, where an edge $\{i,j\}$ is proper if $i$ is the parent of $j$ and  $i$ is less than all the descendants of $j$. Despite the similarity in definitions, the distribution of proper edges in trees is very different from that of proper vertices.

\newsubsection

In Section 2, we prove some lemmas relating proper vertices of forests  to hook length formulas. Sections 3 and 4 contain background information on exponential generating functions and the enumeration of trees and forests. Some properties of the polynomials $P_n(a,b,c)$, including a differential equation satisfied by its generating function, are proved in Section 5.
The next three sections contain our main results on the enumeration by proper vertices of forests of rooted trees, $k$-ary trees, and forests of ordered trees, with corresponding hook length formulas. In Section 9 we discuss another interpretation of $P_n(a,b,c)$, due to 
E\u gecio\u glu and Remmel \cite{eg-rem}, in terms of ascents and descents in forests, and in Section 10 we give a parking function interpretation to a specialization of $P_n(a,b,c)$.

\section{Hook length formulas}
\label{s-hook}
Throughout this paper, all trees are rooted, and all forests are forests of rooted trees. An \emph{ordered tree} is a tree in which the children of each vertex are linearly ordered. An \emph{ordered forest} is a forest of ordered trees in which the trees are also linearly ordered. We write $V(F)$ for the set of vertices of the forest $F$.

A vertex $v$  in a tree is a \emph{descendant} of a vertex $u$ if $u$ lies on the unique path from the root to $v$. If $v$ is a descendant of $u$, it is a \emph{proper descendant} if it is not equal to~$u$.

Let $F$ be a forest with $n$ vertices. Then a \emph{labeling} of $F$ is a bijection from $V(F)$ to~$[n]$. A vertex of $F$ is \emph{proper} with respect to a labeling of $F$ if its label is less than the labels of all its proper descendants. The \emph{hook length} $h(v)$ of a vertex $v$ in a forest is the number of descendants of $v$, including $v$.

We will consider various kinds of trees and forests in this paper, and for each of them we will have a ``hook length formula" something like \eqref{e-P2}. These hook length formulas follow from our results on counting forests by proper vertices and the  following lemmas, which relates hook lengths to proper vertices.

\begin{lemma}
\label{l-1}
Let $F$ be a forest with  $n$ vertices, and let $S$ be a subset of $V(F)$. The the number of labelings of $F$ in which every vertex in $S$ is proper is 
\begin{equation*}
\frac{n!}{\prod_{v\in S}h(v)}.
\end{equation*}
\end{lemma}

\begin{proof}
If all labelings of $F$ are equally likely, then the probability that a vertex $v$ is proper is $1/h(v)$. It can be shown that these probabilities 
are independent, and this implies the lemma.
\end{proof}

\begin{lemma}
\label{l-hook1}
Let $F$ be a forest with $n$ vertices. Then
\begin{equation}
\label{e-hook1}
\sum_{L} (1+\alpha)^{\prop L} =n!\!\!\prod_{v\in V(F)}
\left(1+\frac {\alpha}{h(v)}\right),
\end{equation}
where $L$ runs over all labelings of $F$.
\end{lemma}

\begin{proof}
The left side of \eqref{e-hook1} is equal to 
\begin{math}
\sum_{L,S} \alpha^{|S|}
\end{math}
where $L$ runs over all labelings of $F$ and $S$ runs over all subsets of the set of vertices of $F$ that are proper with respect to $L$.
The right side of \eqref{e-hook1} is equal to
\begin{equation*}
\sum_{S\subseteq V(F)}\frac{n!}{\prod_{v\in S} h(v)}\alpha^{|S|}.
\end{equation*}
The result then follows from Lemma \ref{l-1}.
\end{proof}


Let us say that two labelings of a forest $F$ are \emph{isomorphic} if there is an automorphism of $F$ that takes one labeling to the other. (An automorphism of an ordered forest must preserve the order structure, so ordered forests have no nontrivial automorphisms.) Let $F$ be a forest on $n$ vertices with automorphism group $G$. Then the $n!$ labelings of $F$ fall into $n!/|G|$ isomorphism classes. It is clear that isomorphic labelings have the same number of proper vertices, so we obtain the following lemma by dividing Lemma \ref{l-hook1} by $|G|$.

\begin{lemma}
\label{l-reps}
Let $F$ be a forest with $n$ vertices and let $\mathcal L$ be a set of representatives of isomorphism classes of labelings of $F$. Then 
\begin{equation}
\label{e-iso}
\sum_{L\in \mathcal L} (1+\alpha)^{\prop L}=|\mathcal L|\!\! \prod_{v\in V(F)}\left(1+\frac\alpha{h(v)}\right).
\end{equation}
\end{lemma}
Let us now define a \emph{labeled forest} to be a forest whose vertex set is $[n]$ for some $n$. Thus replacing each vertex of a labeling of an arbitrary forest gives a labeled forest, and two labelings give the same labeled forest if and only if the labelings are isomorphic.

Then \eqref{e-iso}, together with linearity, gives the following result.

\begin{lemma}
\label{l-hook3}
Let $\mathcal F$ be a class of labeled forests with vertex set $[n]$ in which membership depends only on isomorphism class. Then 
\begin{equation*}
\sum_{F\in \mathcal F}(1+\alpha)^{\prop F}=\sum_{F\in \mathcal F}
\prod_{v\in [n]}\left(1+\frac\alpha{h(v)}\right).
\end{equation*}
\end{lemma}

For ordered forests, which have no nontrivial automorphisms, we can simplify this result, since each isomorphism class of ordered forests on $n$ vertices corresponds to $n!$ labeled ordered forests. An \emph{unlabeled forest} may be formally defined as an isomorphism class of forests, though we think of an unlabeled forest as a forest in which the identity of the vertices is irrelevant. Then the following lemma follows from either Lemma \ref{l-reps} or Lemma \ref{l-hook3}:

\begin{lemma}
\label{l-hook}
Let $\mathcal A$ be a set of unlabeled ordered forests with $n$ vertices,
and let $\mathcal F$ be the corresponding set of labeled ordered forests. Then 
\begin{equation}
\label{e-hook}
\sum_{F\in \mathcal F}(1+\alpha) ^{\prop F} 
= n! \sum_{F\in \mathcal A}\prod_{v\in V(F)}
  \left(1+\frac\alpha{h(v)}\right),
\end{equation}
where the product is over all vertices $v$ of $F$.
\end{lemma}

In later sections, it will be convenient to use the term ``labeled forest" to refer more generally to any forest whose vertices are integers.

\section{Exponential generating functions}
\label{s-egf}

We assume that the reader is familiar with the combinatorial interpretation of basic operations on exponential generating functions: addition, multiplication, composition, and differentiation, as described, for example in \cite[Chapter 5]{ec2} or \cite{bll}. 
Rather than adapting the formal language of \cite{bll}, we take a more informal approach to the connection between exponential generating functions and the objects that they count. Thus in deriving functional or differential equations for generating functions we will describe the decompositions that lead to these equations, but omit the straightforward details of how the equations follow from the decompositions. All of our exponential generating functions will be in the variable~$x$, so we will usually omit it, and write, for example $A$ instead of $A(x)$. Since nearly all our generating functions are exponential, we will usually call them simply ``generating functions"; when we use ordinary generating functions we will note this explicitly.

There are two kinds of decompositions of labeled objects that are especially important: \emph{unordered} and \emph{ordered} decompositions. An \emph{unordered} decomposition is a decomposition of a labeled object into an unordered collection of objects. An important example is the decomposition of a forest  into its component trees. When there is an unordered decomposition of objects counted by a generating function $f$ into objects counted by $g$, the generating functions are related by the ``exponential formula" $f=e^g$. Decompositions may depend on a total ordering of the labels. For example, we can decompose a permutation (considered as a linear arrangement) by cutting it before each left-right minimum, yielding a set of permutations each starting with its smallest element. (Thus the permutation $4 7 5 3 1 2 6$ is decomposed into $\{4 7 5, 3, 1 2 6\}$.) The permutation can be reconstructed from its constituents by arranging them in decreasing order of their first elements.
The generating function identity in this example is
\begin{equation}
\label{e-explog}
\egf{n!} 
=\exp\biggl(\sumo n (n-1)! \dpow xn\biggr).
\end{equation}

An \emph{ordered decomposition} is a decomposition into a sequence of objects. If objects counted by $f$ have an ordered decomposition into objects counted by $h$ then $f$ and $h$ are related by $f=1/(1-h)$.  The simplest example is the decomposition of a permutation, considered as a linear arrangement, into the sequence of its entries; here $g=x$ and $f=1/(1-x)$.
Another important example is the decomposition of an ordered forest into its constituent trees. 

If an ordered decomposition for a class of labeled objects, then an unordered decomposition must also exist: if $f=1/(1-h)$ then $f = \exp(\log (1-h)^{-1})$ and we can interpret this by replacing $x$ with $h$ in the combinatorial interpretation given above for \eqref{e-explog}. We will see an example of this kind of decomposition in Section \ref{s-of}.

\section{Trees and forests}
\label{s-tf}
We introduce here the types of trees and forests that we will count by various parameters in later sections. Most of this material is well known, but it is convenient to clarify the definitions and state the known  formulas that we will generalize later. 

Let $A$ be the generating function for labeled forests.  Then since a tree consists of a root together with the forest of proper descendants of the root, the generating function for trees is $xA$, and since a forest is a set of trees, the generating function for forests is $e^{xA}$. Thus $A$ satisfies the functional equation
\begin{equation}
\label{e-forests}
A= e^{xA}.
\end{equation}
It is well known that the solution of $\eqref{e-forests}$ is 
\[A = \egf{(n+1)^{n-1}},\]
and that more generally,
\begin{equation}
\label{e-c2}
A^c = \egf{c(n+c)^{n-1}}.
\end{equation}
Equivalent formulas can be found in \cite[Sections 1.2.5 and 1.2.7]
{g-j},  \cite[Proposition 5.4.4]{ec2} and \cite[p.~97]{riordan}.
Since $A^c = e^{c\log A} = e^{c xA}$, the coefficient of $c^i$ in $A^c$ is the generating function for forests of $i$ trees.



\begin{figure}[htbp] 
   \centering
   \includegraphics[width=3.5in]{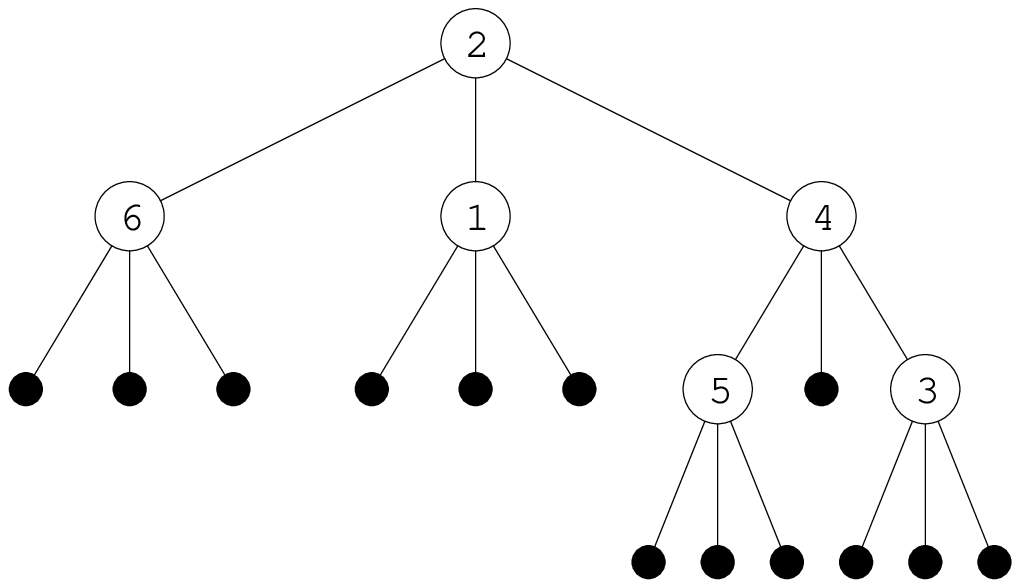} 
   \caption{A labeled ternary tree}
   \label{f-2}
\end{figure}

We say that a vertex of a tree has \emph{degree $d$} if it has $d$ children. 
A $k$-ary tree is an ordered tree in which every vertex has degree $k$ or 0. The ordinary generating function $D$ for unlabeled $k$-ary trees in which the coefficient of $x^n$ is the number of $k$-ary trees with $n$ vertices of degree $k$ (and thus $1+(k-1)n$ vertices of degree 0) satisfies
\begin{equation}
\label{e-k0}
D = 1 +xD^k,
\end{equation}
and it is well known that\begin{equation}
\label{e-k1}
D^c = \sumz n \frac{c}{kn+c} \binom{kn+c}n x^n.
\end{equation}
See, for example, \cite[Exercise 2.7.1]{g-j}, \cite[Proposition 6.2.2]{ec2}, and \cite[pp.~200--201]{gkp}. In computing hook lengths of $k$-ary trees, we will ignore vertices of degree 0. Thus the hook length of the root of the ternary tree in Figure \ref{f-2} is 6.

For our purposes a labeled $k$-ary tree is a $k$-ary tree in which only the vertices of degree $k$ are labeled.
Figure~\ref{f-2} shows a labeled ternary tree.
The exponential generating function for labeled $k$-ary trees is the same, as a power series, as the ordinary generating function for unlabeled $k$-ary trees, but as an exponential generating function we would rewrite 
\eqref{e-k1} as 
\[
D^c = \egf{c (kn+c-1)(kn+c-2)\cdots ((k-1)n+c+1) }
\]

The generating function $E$ for ordered trees (by the total number of vertices) satisfies 
$E = x/(1-E).$
As in the case of $k$-ary trees this formula holds for both the ordinary generating function, which counts unlabeled ordered trees, and for the exponential generating function, which counts labeled ordered trees. The generating function $F$ for ordered forests (in either the unlabeled or labeled version) satisfies 
$F=1/(1-E)$ and $E=xF$, and thus $F=1/(1-xF)$.

The formula $F=1/(1-xF)$ is easily transformed into $F=1+xF^2$, so, as is well known, the generating function for ordered forests is the same as that for binary trees. In order to generalize this to $k$-ary trees, we introduce \emph{$k$-colored ordered forests} with generating function $F$ satisfying
\begin{equation}
\label{e-kcolor}
F=\frac 1{1-xF^k}.
\end{equation}
A $k$-colored ordered forest is a forest of ordered trees in which the edges are colored in colors 1, 2, \dots,  $k$ with the property that  among the children of any vertex, those of color 1 come first, then those of color 2, and so on. (But not every color must appear.) 
\begin{figure}[htbp] 
   \centering
   \includegraphics[width=3in]{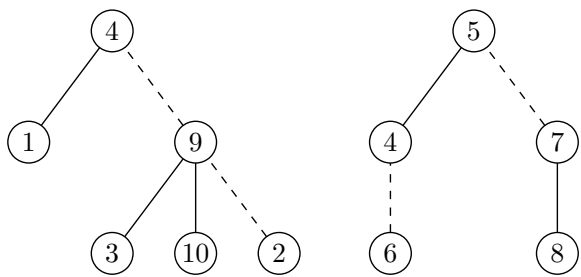} 
   \caption{A 2-colored forest}
   \label{f-2color}
\end{figure}
For example, Figure \ref{f-2color} shows a 2-colored ordered forest in which color 1 is represented by solid lines  and  color 2 by dotted lines.
Since \eqref{e-kcolor} is equivalent to $F=1+xF^{k+1}$, $k$-colored forests are equinumerous with $(k+1)$-ary trees, and the well-known bijection from forests of ordered trees  to binary trees
(see, e.g., Stanton and White \cite[p.~61]{sw}) generalizes easily to this situation. 

\section{The generating function for $P_n(a,b,c)$}

Our main objects of study are the polynomials $P_n(a,b,c)$ defined by 
\begin{equation}
\label{e-1}
P_n(a,b,c)= c\prod_{i=1}^{n-1}(ia+(n-i)b +c).
\end{equation}
Note that $P_n(a,b,c)$ is symmetric in $a$ and $b$. It is also homogeneous of degree of $n$ in $a$, $b$, and $c$, so with no loss of information we could set any of the variables to one (or to any other nonzero value). 

We denote by $Q_n(a,b)$ the coefficient of $c$ in $P_n(a,b,c)$, so  
\begin{equation*}
Q_n(a,b) = c^{-1}P_n(a,b,c)\Big|_{c=0} = 
       \prod_{i=1}^{n-1}(ia+(n-i)b).
\end{equation*}

The polynomials $P_n(a,b,c)$ generalize the counting formulas of the previous section, since
\[c(n+c)^{n-1}=P_n(1,1,c)\]
and 
\begin{equation}
\label{e-P3}
 \frac{c}{kn+c} \binom{kn+c}n  n! = P_n(k-1,k,c).
\end{equation}
Since the polynomial \eqref{e-P3} has two free parameters, by rescaling we can obtain $P_n(a,b,c)$ and we find that for $a\ne b$, we have 
\begin{equation}
\label{e-P4}
P_n(a,b,c)/n! = (b-a)^n\frac{ \bar c}{\bar b n + \bar c}
\binom{\bar b n +\bar c}{n},
\end{equation}
where $\bar b = b/(b-a)$ and $\bar c=c/(b-a)$.

It follows from \eqref{e-P4}, \eqref{e-k0} and \eqref{e-k1} that 
\[
\egf{P_n(a,b,c)}  = G^{c/(b-a)},
\]
where $G$ satisfies 
\begin{equation}
\label{e-Pgf}
G = 1 + (b-a)x G^{b/(b-a)}.
\end{equation}
It does not seem to be easy to obtain a combinatorial interpretation for $P_n(a,b,c)$ directly from \eqref{e-Pgf}. Instead, we derive a differential equation whose combinatorial meaning is clearer. Let $H=G^{1/(b-a)}$, so that 
\begin{equation}
\label{e-PH}
\egf{P_n(a,b,c)} = H^c.
\end{equation}
We claim that $H$ satisfies the differential equation
\begin{equation}
\label{e-H3}
H' = H^{a+1}+bxH^a H'.
\end{equation}
To see this, \eqref{e-Pgf} gives

\begin{equation}
\label{e-H1}
H^{b-a} = 1+(b-a)xH^b.
\end{equation}
Differentiating with respect to $x$ and dividing both sides by $(b-a)H^{b-a-1}$ gives \eqref{e-H3}.
(Although this derivation is not valid for $a=b$, a limit argument shows that \eqref{e-H3} does hold in this case.)
Conversely, \eqref{e-H3}, together with the initial value $H(0)=1$, implies \eqref{e-H1}, and thus \eqref{e-PH}.

We note that \eqref{e-H1} can be written in the more symmetric and combinatorially suggestive form 
\begin{equation}
\label{e-Hsym}
H^b(1+axH^{a})=H^a(1+bxH^{b});
\end{equation}
however, we have not found a combinatorial interpretation to 
\eqref{e-Hsym}. 

We now give a completely self-contained proof that \eqref{e-H3} implies \eqref{e-PH}. For our applications, it is convenient to 
add another redundant variable to the differential equation.



\begin{theorem} 
\label{t-1'}
Let $g=g(x)$ be the solution of the differential equation
\begin{equation}
\label{e-1'}
g'=ug^{b+1} +axg^b g'
\end{equation}
with $g(0)=1$, where the prime denotes the derivative with respect to $x$. Then 
\begin{equation}
\label{e-gc}
g^c = \egf{P_n(a,bu,cu)}.
\end{equation}
and
\begin{equation}
\label{e-Q}
\log g = \sumo n uQ_n(a,bu)\dpow xn.
\end{equation}
\end{theorem}

\begin{proof}
\renewcommand{\P}{\bar P}

Let  $\P_n(c)$ be the coefficient of $x^n\!/n!$ in $g^c$, where the dependence of $\P_n(c)$ on $u$, $a$, and $b$ is implicit. 
From \eqref{e-1'} we obtain
\begin{equation}
\label{e-2'}
(g^c)' = c u g^{b+c} + acxg^{b+c-1} g' = cug^{b+c} 
       + \frac{acx}{b+c}(g^{b+c})'
\end{equation}
Equating coefficients of $x^n\!/n!$ in \eqref{e-2'} gives
\begin{align*}
\P_{n+1}(c) &= cu\P_n(b+c)+\frac{nac}{b+c}\P_n(b+c)\\
    &=\frac{c(na+bu+cu)}{b+c}\P_n(b+c).
\end{align*}    
Since  $\P_0(c)=1$, we get 
$$\P_n(c)=cu(a+(n-1)bu+cu)(2a+(n-2)bu+cu)\cdots((n-1)a+bu+cu).$$
The formula for $\log g$ follows by taking the coefficient of $c$ in $g^c$.
\end{proof}

\section{Counting forests by proper vertices}

As noted earlier, $P_n(1,1,1) = (n+1)^{n-1}$ is the number of forests on $[n]$. Thus we might hope that $P_n(a,b,c)$ counts these forests according to some parameters. Theorem~\ref{t-2} shows that this is the case.

For a forest $F$, let $\prop F$ be the number of proper vertices and let $\tree F$ be the number of trees in $F$. We give two proofs of the following result, one using the differential equation for the generating function for $P_n(a,b,c)$ and one which is more combinatorial, though not completely bijective.

\begin{theorem}
\label{t-2}
We have 
\begin{equation}
\label{e-T1b}
\sum_F a^{n-\prop F} b^{\prop F -\tree F}c^{\tree F}=P_n(a,b,c),
\end{equation}
where the sum is over all forests $F$ on $[n]$.

\end{theorem}

\begin{proof}[First proof]
Let $I$ be the generating function for forests by the number of proper vertices; i.e., 
\begin{equation*}
I=\sum_{n=0}^\infty \sum_F u^{\prop F} \dpow xn
\end{equation*}
where $F$ runs over all forests on $[n]$.




Let $J$ be the generating function for trees weighted by proper vertices. Then $I=e^J$.  Next we will show that $J' = uI +xI'$.
Let $T$ be a tree on $\{0,1,\dots, n\}$, interpreted as a labeled object with label set $[n]$. Then $T$ consists of a root together with the forest $F$ obtained by deleting the root. If the root of $T$ is 0 then $T$ has one more proper vertex than $F$. Thus the contribution to $J'$ from trees in which the root is 0 is $uI$. If the root is not 0, then $T$ has an improper root with a label in $[n]$, and 0 appears in $F$. Thus the contribution to $J'$ from trees in which the  root is not 0 is $xI'$.

Then \[I' = (e^J)' = e^J J' = I(uI+xI') = uI^2 + xII',\]
so by Theorem \ref{t-1'}, 
\[
I= \egf{P_n(1,u,u)},
\]
as shown by Seo \cite[Corollary 8]{seo}.

The number of proper vertices of a forest is the same as the number of proper vertices in all its constituent trees, so 
 $I^c = \exp( cJ)$ counts forests  by proper vertices in which every tree is weighted by $c$, which gives 
\begin{equation}
\label{e-T1a}
\sum_F u^{\prop F}c^{\tree F} = P_n(1,u,cu)
\end{equation}
Using the fact that $P_n(a,b,c)$ is a homogeneous polynomial of degree $n$ in $a$, $b$, and $c$, we can rescale \eqref{e-T1a} to get \eqref{e-T1b}.
\end{proof}

We can interpret \eqref{e-T1b} as counting forests on $[n]$ in which each vertex is assigned a weight. We will call a vertex 
\emph{weakly proper} if it is proper, but not the smallest vertex in its tree and \emph{minimal} if it is the smallest vertex in its tree. Then \eqref{e-T1b} counts forests in which every improper vertex is weighted $a$, every weakly proper vertex is weighted $b$, and every minimal vertex is weighted $c$.

\begin{proof}[Second proof of Theorem \ref{t-2}]
Let $p_n(a,b,c)$ count forests $F$ on $[n]$ in which every improper
vertex is weighted $a$,
every weakly proper vertex is weighted $b$, and every minimal vertex in a
component of $F$ is weighted $c$.
Then $c^{-1}p_n(a,b,c)$ counts forests with the same weight on each
vertex
except that the smallest vertex of $F$ is weighted by $1$.
Suppose that a forest $F$ on $[n]$ is weighted in this way.
If $1$ is a root of $F$, delete the vertex $1$ of $F$.
Then we obtain a forest on $\{2,3,\ldots,n\}$ in which the minimal
vertices
are weighted either $b$ or~$c$. (See (a) of Figure~\ref{f-1}.)
\begin{figure}[htbp]    \centering
   \includegraphics[width=5.5in]{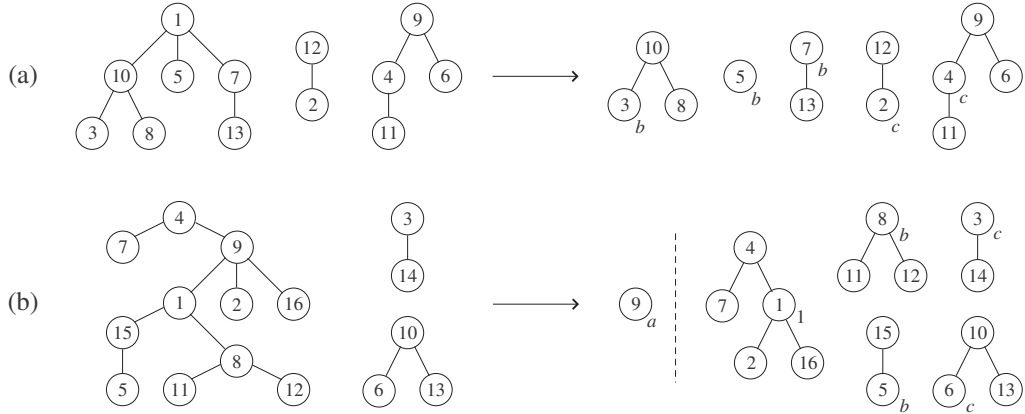}
   \caption{A weight preserving decomposition of forests}
   \label{f-1}
\end{figure}
Thus these forests are counted by $p_{n-1}(a,b,b+c)$.
If $1$ is not a root of $F$, interchange the label $1$ and its parent $i$
and delete the vertex $i$.
Then we have a forest on $[n]\setminus\{i\}$ such that the smallest
vertex $1$ is weighted $1$ and
the other minimal vertices are weighted either $b$ or $c$.
Moreover, there are $n-1$ possibilities for the parent $i$ of~$1$  and it is weighted $a$.
Thus these forests are counted by
$(n-1)\,a\cdot(b+c)^{-1}p_{n-1}(a,b,b+c)$. (See (b) of Figure~\ref{f-1}.)
Combining two cases yields
$$
c^{-1}p_n(a,b,c)=p_{n-1}(a,b,b+c)+(n-1)\,a\,(b+c)^{-1}p_{n-1}(a,b,b+c),
$$
with the initial condition $c^{-1}p_1(a,b,c)=1$. It follows easily that
\[p_n(a,b,c)=P_n(a,b,c).\]

\end{proof}


The symmetry $P_n(a,b,c)=P_n(b,a,c)$ implies that among forests of $j$ trees, the number with $i$ improper vertices is equal to the number with $i$ weakly proper vertices. (It is clear that the case $j=1$ of this symmetry implies the general case.) We do not have a bijective proof of this symmetry.

 It is instructive to look at the case $i=0$ of this symmetry, which can be explained by known bijections. By the homogeneity of the polynomials $P_n(a,b,c)$, the symmetry is equivalent to the identity $P_n(0,1,c)=P_n(1,0,c)=c(c+1)\cdots(c+n-1)$, in which the coefficients are unsigned Stirling numbers of the first kind. First, $P_n(0,1,c)$ counts forests on $[n]$ with $n$ proper vertices. These are increasing forests, and there is a well known bijection \cite[p.~25]{ec1}
that takes an increasing forest of $j$ trees to a permutation with $j$ left-right minima. 
Next, $P_n(1,0,c)$ counts forests in which the number of proper vertices equals the number of trees. Since every tree has at least one leaf, and every leaf is proper, in these forests every tree has a single leaf, which is the least vertex of the tree. If we read the trees from leaf to root, we get a set of $j$ permutations, each of which starts with its least element, and as described in Section \ref{s-egf} this set of permutations can be arranged to form a permutation with $j$ left-right minima.

A different symmetry of forests that has a very similar special case is described in \cite{gessel}. (See also \cite{kalikow}.)

Since the proper vertices of a forest are the same as the proper vertices of its component trees, \eqref{e-Q} gives a corresponding result for trees;
\begin{cor}
We have
\begin{equation*}
\sum_T a^{n-\prop T}b^{\prop T -1}= Q_n(a,b),
\end{equation*}
where $T$ runs over all trees on $[n]$.\qed
\end{cor}


Applying Lemma \ref{l-hook3} gives a hook-length formula for forests:
\begin{cor}
\label{c-forest-hook}
\begin{equation*}
\sum_{F} c^{\tree F}\prod_{v\in [n]}
  \left(1+\frac\alpha{h(v)}\right)
 = P_n(1, 1+\alpha, c(1+\alpha)),
 \end{equation*}
 where $F$ runs over all forests on $[n]$.\qed
 \end{cor}
 In particular, taking $\alpha=1$ gives
 \begin{equation*}
\frac1{n!}\sum_{F\in \mathcal F} c^{\tree F}\prod_{v\in [n]}
  \left(1+\frac1{h(v)}\right)
 = \frac1{n!}\,P_n(1, 2, 2c) = \frac c{n+c} \binom{2n+2c}{n},
\end{equation*}
which for $c=1$ is the Catalan number $C_{n+1}$.

\section{$k$-ary trees}

In this section we give an analogue of Theorem \ref{t-2} for $k$-ary trees. 
We first count $k$-ary trees by proper vertices, obtaining a result of
Seo \cite[Corollary 7]{seo}.
\begin{theorem}
\label{t-k-ary}
We have
\begin{equation}
\sum_T u^{\prop T}= P_n(k,(k-1)u, u),
\end{equation}
where the sum is over all $k$-ary trees $T$ on $[n]$.
\end{theorem}

\begin{proof}
Let $J$ be the generating function for $\sum_T u^{\prop T}$ as in the statement of the theorem. By reasoning similar to that used in the proof of Theorem \ref{t-2}, we will obtain the differential equation 
\begin{equation}
\label{e-k-ary}
J' = u J^k + kxJ^{k-1}J'.
\end{equation}
We consider $k$-ary trees on the vertex set $\{0,1,\dots, n\}$ and look at the location of vertex~0. The term $uJ^k$ in \eqref{e-k-ary} corresponds to the case in which 0  is the root of the tree. If 0 is not the root, the tree can be decomposed into a root $R$ (the factor $x$), the subtree containing~0 rooted at a child of $R$ (the factor $J'$) and a $(k-1)$-tuple of the other subtrees rooted at the children of $R$ (the factor $J^{k-1}$). The factor $k$ is the number of positions  in which the subtree containing 0 can 
lie in relation to the other subtrees.

So by Theorem \ref{t-1'},
\begin{equation}
\label{e-ktree}
J = \egf {P_n(k, (k-1)u, u)}.
\end{equation}
\end{proof}



Applying Lemma \ref{l-hook1}, gives the following hook-length formula for $k$-ary trees, which is equivalent to Corollary 3.11 of Du and Liu \cite{du-liu}. The case $k=2$, $\alpha = 1$ is equivalent to Postnikov's formula \eqref{e-P2}.

\begin{cor}
\label{c-k-ary-hook}
We have
\begin{equation*}
n!\sum_T \prod_{v\in V(T)} \left(1+\frac {\alpha}{h(v)}\right)
 = P_n(k,(k-1)(1+\alpha ), 1+\alpha)
\end{equation*}
where $T$ runs over all unlabeled $k$-ary trees on $n$ vertices.
\end{cor}

\begin{figure}[htbp] 
   \centering
   \includegraphics[width=3in]{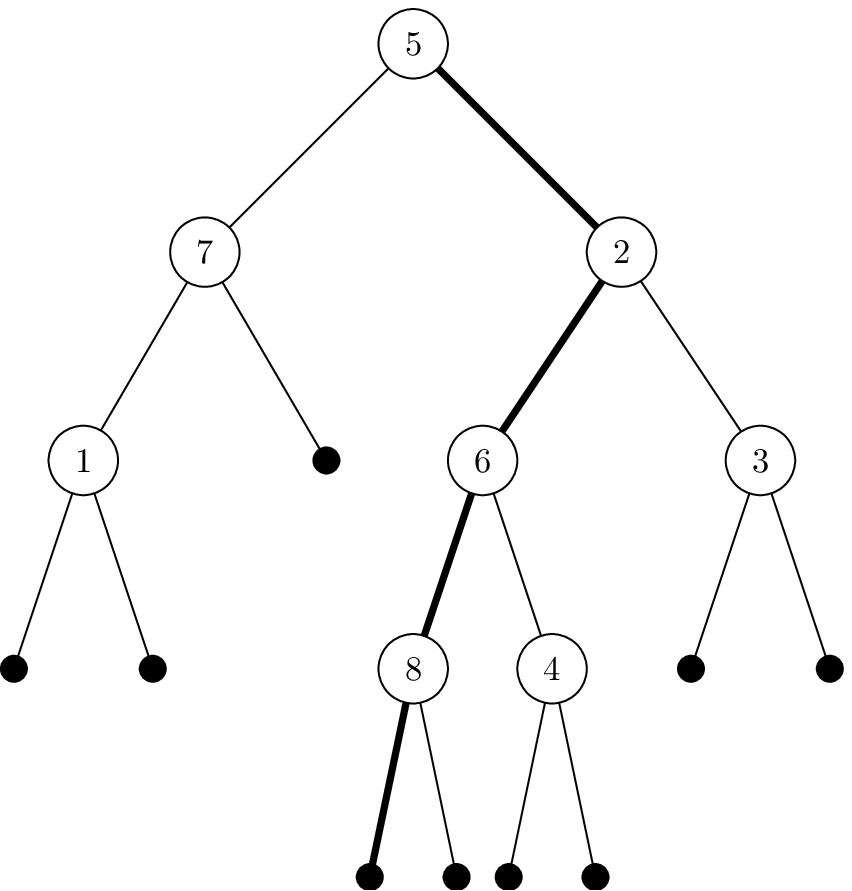} 
   \caption{A binary tree and its \dpath}
   \label{f-dpath}
\end{figure}

\begin{figure}[htbp] 
   \centering
   \includegraphics[width=3in]{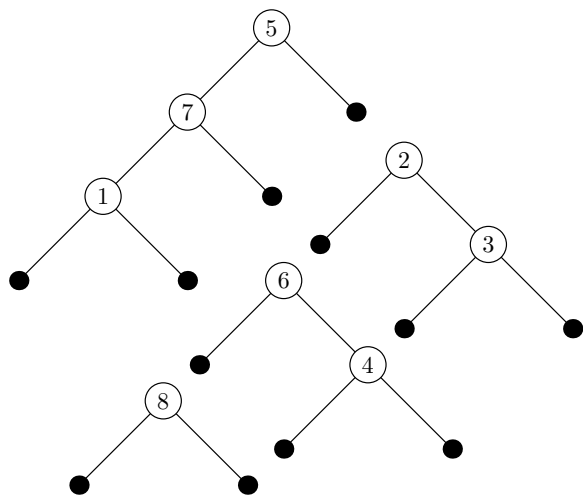} 
   \caption{Decomposition of the tree of Figure \ref{f-dpath}}
   \label{f-decomposed}
\end{figure}

We can generalize Theorem \ref{t-k-ary} to include a parameter that corresponds to the number of trees in a forest, and thereby find a result analogous to Theorem \ref{t-2}. In order to do this, we must find an unordered decomposition of $k$-ary trees that is compatible with proper vertices; i.e., no improper vertex of a $k$-ary tree becomes proper when restricted to its component. To define this decomposition we first specify a particular path in a $k$-ary tree $T$ from the root to a leaf. The rule is that we always go from a labeled vertex $v$ to the first child of $v$ unless $v$ is improper and the smallest descendant of $v$ is also a descendant of the first child of $v$. In this case we go to the second child of $v$. When we reach a leaf, we stop. We call this path the \emph{\dpath} of $T$. For example, Figure \ref{f-dpath} shows a binary tree with its 
\dpath.
The unordered decomposition of the $k$-ary $T$ is obtained by removing each edge on its \dpath\ (except the last edge ending at a leaf) and replacing each lost child with a leaf. We
denote by $\comp T$ the number  components in the decomposition of $T$, which is the same as the number of labeled vertices in the \dpath. 
 The unordered decomposition of the tree of Figure \ref{f-dpath} is shown in Figure \ref{f-decomposed}.

\begin{theorem}
\label{t-kcomps}
For $k>1$ we have
\begin{equation}
\label{e-kary0}
\sum_{T} v^{n-\prop T}u^{\prop T -\comp T}w^{\comp T}
=P_n(kv,(k-1)u,  w).
\end{equation}
where the sum is over all $k$-ary trees $T$ on $[n]$.
\end{theorem}
\begin{proof} 
We first prove
\begin{equation}
\label{e-kary1}
\sum_T u^{\prop T}c^{\comp T}= P_n(k,(k-1)u, cu).
\end{equation}

By \eqref{e-gc} and Theorem \ref{t-k-ary},  to prove \eqref{e-kary1} it is sufficient to show that the generating function for the left side is $J^c$, where $J$ is as in the proof of Theorem  \ref{t-k-ary}. 
This follows from a verification that the unordered decomposition for 
a $k$-ary tree has several properties. Let the decomposition of a $k$-ary  tree $T$ be $\{T_1, \dots, T_m\}$. We must show that $\comp T_i= 1$ for each $i$, that $T$ can be recovered from $\{T_1, \dots, T_m\}$, and that $\prop T = \sum_i \prop T_i$. Moreover we must also show that any such set of trees $\{T_1, \dots, T_m\}$ with disjoint labels comes from the decomposition of some $k$-ary tree.

The components $T_i$ are characterized by the  property that either the first child of $T_i$ is a leaf, or the root of $T_i$ is improper, its smallest descendant is a descendant of its first child, and its second child is a leaf. This implies that $\comp T_i =1$. We can reconstruct the order of the components since 
the first component $T_1$  will contain the smallest vertex of $T$, the next component $T_2$ will contain the smallest vertex of $T$ not in $T_1$, and so on. To reconstruct $T$, let $r_i$ denote the  root of tree $T_i$. Then in $T$, $r_{i+1}$ must have been either the first or second child of $r_{i}$: $r_{i+1}$ was the first child of $r_i$ in $T$ if and only if the first child of $T_i$ is a leaf. The decomposition preserves proper and improper vertices since the smallest descendant of any improper vertex $v$ on the \dpath\ remains a descendant of $v$ in its component (and descendants of vertices not on the \dpath\ are not changed). The last assertion is clear.

The homogeneity of $P_n(a,b,c)$, together with \eqref{e-kary1}, implies that 
\begin{equation*}
\sum_T v^{n-\prop T}u^{\prop T}c^{\comp T}= P_n(kv,(k-1)u, cu).
\end{equation*}
Replacing $c$ with $w/u$ gives \eqref{e-kary0}.
\end{proof}

In should be pointed out that Theorem \ref{t-kcomps} is not valid for unary trees. However, there is a different (and simpler) decomposition for unary trees that gives an analogue of Theorem
\ref{t-kcomps} in this case.

In contrast to Corollary \ref{c-forest-hook}, we cannot use Theorem \ref{t-kcomps} to add a parameter for components to Corollary \ref{c-k-ary-hook}, since the number of components of a $k$-ary tree is not an isomorphism invariant.

\section{Ordered forests} 
\label{s-of}

Next we count $k$-colored ordered forests by proper vertices. The case $k=1$ of the next theorem was given by Seo \cite[equation (14)] {seo}.


\begin{theorem}
\label{t-kcolor0}
We have
\begin{equation}
\label{e-kcolor1}
\sum_{F} u^{\prop F}=P_n(k,(k+1)u,u),
\end{equation}
where the sum is over all $k$-colored ordered forests $F$ on $[n]$,
and for $n\ge 1$, 
\begin{equation}
\label{e-kcolor2}
\sum_{T}u^{\prop T} = -P_n(k, (k+1)u, -u),
\end{equation}
where the sum is over all $k$-colored ordered trees $T$ on $[n]$.
\end{theorem}

\begin{proof}
Let $K$ be the generating function for $k$-colored ordered forests, where a forest $F$ is weighted $u^{\prop F}$. Let $L$ be the analogous generating function for $k$-colored ordered trees, so $K=1/(1-L)$.
We claim that 
\begin{equation}
\label{e-KL}
L' = uK^k+kxK^{k-1}K'.
\end{equation}
By similar reasoning to that used before, we see that in \eqref{e-KL} the term $uK^k$ counts $k$-colored ordered trees on $\{0,1,\dots, n\}$ in which 0 is the root, and for each $i$ from 1 to $k$, 
$xK^{k-1}K'$ is the generating function for forests in which 0 is in a forest joined to the root by an edge of color $i$. 
It follows from \eqref{e-KL} and $K=1/(1-L)$ that 
\[K' = uK^{k+2} + kxK^{k+1}K'.\]
So by Theorem \ref{t-1'} we have 
\begin{equation}
\label{e-K}
K=\egf{P_n(k,(k+1)u,u)},
\end{equation}
proving the first formula. 

For the second formula, we note that
\[L=1-K^{-1},\]
and a formula for the coefficients of $K^{-1}$ follows from 
\eqref{e-K} and \eqref{e-gc}.

\end{proof}

Applying Lemma \ref{l-hook} to Theorem \ref{t-kcolor} gives hook length formulas for $k$-colored ordered forests and trees. A result equivalent to the case $k=1$ of \eqref{e-kchook} was given by Du and Liu \cite[Theorem 4.14]{du-liu}.

\begin{cor} We have
\begin{equation}
\label{e-kchook}
n!\sum_F \prod_{v\in V(F)}\left(1+\frac\alpha{h(v)}\right)
  = P_n(k, (k+1)(1+\alpha),1+\alpha),
\end{equation}
where the sum is over all unlabeled $k$-colored ordered forests $F$ on $n$ vertices, and 
\begin{equation*}
n!\sum_T \prod_{v\in V(T)}\left(1+\frac\alpha{h(v)}\right)
  = -P_n(k, (k+1)(1+\alpha),-(1+\alpha)),
\end{equation*}
where the sum is over all unlabeled $k$-colored ordered trees $F$ on $n$ vertices.
\end{cor}

From the ordered decomposition of $k$-colored ordered forests into $k$-colored ordered trees, we get an unordered decomposition as described in Section \ref{s-egf}. One way to do this is to look at the sequence of roots, and cut the forest after each left-right minimum in this sequence.  (Instead of using the root of each tree for comparison, we could use the smallest  or largest element of each tree.) Thus we decompose each $k$-colored ordered forest into components which are $k$-colored ordered forests in which the first root is the smallest. 
For a $k$-colored ordered forest $F$, let $\comp F$ be the number of components in this sense; that is, the number of left-right minima in the sequence of roots.
Then since 
\begin{equation*}
K^c=\egf{P_n(k,(k+1)u,cu)},
\end{equation*}
we have
\begin{equation*}
P_n(k,(k+1)u,cu) = \sum_F u^{\prop F}c^{\comp F}.
\end{equation*}

By the same reasoning as in Theorem \ref{t-kcomps}, this leads to the analogous formula for $k$-colored ordered trees:
\begin{theorem}
\label{t-kcolor}
We have
\begin{equation*}
\sum_F u^{n-\prop F} v^{\prop F -\comp F} w^{\comp F}
   =P_n(ku, (k+1)v, w),
\end{equation*}
where the sum is over all $k$-colored ordered forests on $[n]$.
\end{theorem}

As noted in Section \ref{s-tf}, there is a simple bijection from $k$-colored ordered forests to $(k+1)$-ary trees, but this bijection does not 
preserve proper vertices. However, we do have a relation between proper vertices of these forests and trees:

\begin{cor}
\label{c-ksym}
We have
\begin{equation*}
\sum_{T} v^{n-\prop T}u^{\prop T -\comp T}w^{\comp T}
= \sum_F u^{n-\prop F} v^{\prop F -\comp F} w^{\comp F}
\end{equation*}
where $T$ runs over $(k+1)$-ary trees on $[n]$ and $F$ runs over $k$-colored ordered forests on $[n]$.
\end{cor}
\begin{proof}
By Theorem \ref{t-kcolor}, the sum on $F$ is equal to $P_n(ku, (k+1)v, w)$.
In \eqref{e-kary0} we replace $k$ by $k+1$, switch $u$ and $v$, and 
use  $P_n(a,b,c)=P_n(b,a,c)$ to show that the sum on $T$ is also equal to 
$P_n(ku, (k+1)v, w)$.
\end{proof}

\section{Descents in trees}
In this section, we describe another forest interpretation, due to 
E\u gecio\u glu and Remmel \cite{eg-rem}, for $P_n(a,b,c)$.
If vertex $i$ is the parent of vertex $j$ in a tree, we call the edge $\{i,j\}$ an \emph{ascent} if $i<j$ and a \emph{descent} if $i>j$, and for a forest $F$ we write $\asc F$ for the number of ascents of $F$ and $\des F$ for the number of descents of $F$. (A different notion of descents in trees was studied in \cite{gessel}.) We use the differential equation approach to count forests by ascents and descents; 
E\u gecio\u glu and Remmel \cite{eg-rem} gave a bijective proof of a more general result. 

\begin{theorem}
We have
\begin{equation}
\label{e-des}
\sum_{F} a^{\des F} b^{\asc F} c^{\tree F}=P_n(a,b,c),
\end{equation}
where the sum is over all forests $F$ on $[n]$.
\end{theorem}

\begin{proof}
Since the ascents and descents of a forest are the same as those of its constituent trees, it is sufficient to prove the case $c=1$ of $\eqref{e-des}$. 
Let $M$ be the generating function for the polynomials $\sum_F a^{\des F} b^{\asc F}$. We shall show that $M$ satisfies the differential equation

We shall show that $M$ satisfies the differential equation
\[M' = M^{b+1} + axM^b M', \]
which by Theorem \ref{t-1'} implies the result.

Let $N = \log M$, so that $M$ counts trees by descents and ascents. Now let $F$ be a forest 
on the vertex set $\{0,1,\dots, n\}$. We decompose $F$ into forests $F_1$ and $F_2$ that together contain all the vertices of $F$ other than 0. We let $F_1$ be the forest of proper descendants of 0, and we let $F_2$ be the forest of all vertices that are not descendants of 0. In order to reconstruct $F$ from $F_1$ and $F_2$ we consider two cases. First, suppose that 0 is a root of $F$. Then $F$ is determined by $F_1$ and $F_2$ and the ascents and descents of $F$ are the same as those of $F_1$ and $F_2$, except that $F$ has one additional ascent for each edge joining 0 to one of its children. The contribution to $M'$ from forests in which 0 is a root is thus $e^{bN}M=M^b M=M^{b+1}$.

If 0 is not a root of $F$, then to reconstruct $F$ from $F_1$ and $F_2$ it is sufficient to know the parent of 0 in $F$; thus $F$ is determined by $F_1$, $F_2$, and the choice of a vertex in $F_2$. All of the ascents and descents of $F_1$ and $F_2$ will be ascents and descents of $F$, and as in the first case, $F$ will have an extra ascent for every child of 0. In addition, $F$ will have an extra descent for the edge joining 0 to its parent. The generating function for forests with a marked vertex is $xM'$, so the contribution to $M'$ from forests in which 0 is not a root is $aM^b \cdot xM'$.
\end{proof}

\section{Parking functions}

A parking function on an $n$-set $K$ is a function $f: K\to [n+1]$ satisfying a certain property, described below. We think of the elements of $K$ as cars and the elements of $[n+1]$ as parking spaces and we think of $f$ as a function that assigns to each car its favorite parking space. The cars arrive in some order at their favorite parking spaces. If space $f(z)$ is not already occupied when car $z$ arrives there, then car $z$ parks in space $f(z)$. More generally, if there is some available space $s$, with $f(z) \le s < n+1$, then car $z$ parks in the first (i.e., smallest) such space. Otherwise, car $z$ does not park at all. If all cars are able to park, then we call $f$ a \emph{parking function}. Although cars are not allowed to park in space $n+1$, we will see that the existence of this space is useful.
It is well known (and easy to see) that whether all cars are able to park or not does not depend on their arrival order. However, we will assume that the cars are  totally ordered and that they arrive in increasing order.

It is useful to represent the parking function $f: K \to [n+1]$ by the ordered partition $(B_1, \dots, B_{n+1})$ of $K$ in which $B_i=f^{-1}(i)$; i.e., $B_i$ is the set of cars that prefer to park in space $i$. 
A necessary and sufficient condition for $f$ to be a parking function is that for each $i$, $|B_1|+\cdots + |B_i|\ge i$.

It is well known that the number of parking functions on $[n]$ is $(n+1)^{n-1}$ \cite[p.~140]{ec2}. So the exponential generating function for parking functions is \[A=\egf{(n+1)^{n-1}}.\]

We  consider a slight generalization of parking functions. Instead of taking the set of parking spaces for $n$ cars to be $[n+1]$, we take it to be $[n+c]$, for some positive integer $c$ (but we keep the number of cars at $n$). Then we call $f$ a \emph{$c$-parking function} if all cars are able to park using the same parking procedure as before, except that cars may now park in spaces 1 to $n+c-1$.
An alternative characterization of $c$-parking functions is that with $B_i = f^{-1}(i)$, we have $|B_1|+\cdots + |B_i|\ge i-c+1$.
It is clear that if all cars park, they will leave exactly $c$ empty spaces, one of which will be the last space, $n+c$. Thus if we represent a $c$-parking function by an ordered partition, then cutting it after every empty parking space decomposes it into a concatenation of $c$ ordinary parking functions.
Therefore the generating function for $c$-parking functions is 
\[A^c = \egf{c(n+c)^{n-1}}.\]
Further properties of $c$-parking functions can be found in Yan \cite{yan}.



If $f$ is a parking function $K\to [n+1]$ then we call a car $z\in K$ \emph{lucky} if it
succeeds in parking in its preferred space. Let $\lucky f$ denote the number of lucky cars associated to~$f$.

\begin{theorem}
\label{t-parking}
We have
\begin{equation*}
\sum_f u^{\lucky f}=P_n(1,u,u),
\end{equation*}
where the sum is over all parking functions $f:[n]\to [n+1]$.
\end{theorem}

\begin{proof}
Let $q_n(u) = \sum_f u^{\lucky f}$, where the sum is over all parking functions 
$f$ on $[n]$  and let
$Q= \sum_{n=0}^\infty q_n(u) x^n/n!$.
We will show that $Q$ satisfies the differential equation 
\begin{equation}
\label{e-parking}
Q' = u Q^2 + xQ' Q,
\end{equation}
which by Theorem \ref{t-1'} implies the result.
To do this we evaluate $Q'$ by counting parking functions $f:[n+1]\to [n+2]$. 
The restriction $f^-$ of such a parking function to $[n]$ is a 2-parking function, which leaves two empty spaces, $e+1$ and $n+2$, for some $e$ with $0\le e\le n$. If $n+1$ is lucky in $f$, then $f(n+1)$ must be $e+1$, and the contribution to $Q'$ from this case is $uQ^2$. If $n+1$ is not lucky in $f$ then $f(n+1)$ is 1, 2, \dots, or~$e$.  In this case $f^-$ decomposes into a parking function on $e$ cars, with one of the first $e$ spaces ``marked," followed by a parking function on $n-e$ cars. The  generating function for parking functions with a marked space is $xQ'$  so the contribution to $Q'$ from this case is  $xQ'Q$. Adding together these two contributions gives \eqref{e-parking}.
\end{proof}

Since the number of lucky cars in a $c$-parking function is the total number of lucky cars among its constituents, the generating function for $c$-parking functions by lucky cars is $Q^c$, and we have the following result:

\begin{cor}
For any positive integer $c$, we have
\begin{equation*}
\sum_f u^{\lucky f}=P_n(1,u,cu),
\end{equation*}
where the sum is over all c-parking functions $f:[n]\to [n+c]$.
\end{cor}

There is an ordered decomposition for parking functions that gives us a result analogous to equation \eqref{e-kcolor2}.
We modify the ordered partition
representation of a parking function $f:K\to [n+1]$, where $|K|=n$,  by removing
the empty block $B_{n+1}$. Let us call the result the
\emph{reduced ordered partition} corresponding to $f$.
Then $f$ is called a \emph{prime parking function} if its
reduced ordered partition cannot be decomposed as a
nontrivial concatenation of reduced ordered partitions
for parking functions; equivalently, if $(B_1,
B_2,\dots, B_n)$ is the reduced ordered partition of
$f$, then $|B_1|+|B_2|+\cdots +|B_i| > i$ for $1\le i <
n$. Then every parking function can be decomposed
uniquely into prime parking functions.

\begin{cor}
For $n\ge1$,
\begin{equation*}
\sum_f u^{\lucky f}=-P_n(1,u,-u),
\end{equation*}
where the sum is over all prime parking functions $f$ on $[n]$.
\end{cor}
\begin{proof}
Let $R$ be the generating for prime parking functions $f$ weighted by $u^{\lucky f}$. Then since the number of lucky cars in a parking function is compatible with the decomposition into prime parking functions, we have
$Q=1/(1-R)$, with $Q$ as in the proof of Theorem \ref{t-parking}.
Thus $R=1-Q^{-1}$, and the result follows from Theorem \ref{t-parking} and 
\eqref{e-gc}.
\end{proof}





%


\end{document}